\documentclass[11pt]{article}%
\usepackage{amsthm}
\usepackage{mathrsfs}

\makeatletter
\def\th@plain{%
  \thm@notefont{}
  \itshape 
}
\def\th@definition{%
  \thm@notefont{}
  \normalfont 
}
\makeatother
\usepackage{amsthm}
\usepackage{mathrsfs}
\usepackage{amsmath}
\usepackage{amsfonts}
\usepackage{latexsym}
\usepackage{amssymb}
\usepackage{graphicx}
\usepackage{tikz}
\usepackage{srcltx,graphicx}
\theoremstyle{plain}
\usepackage{multirow}
\usepackage{soul}
\usepackage{tabularx,ragged2e}
\newcolumntype{C}{>{\Centering\arraybackslash}X} 
\usepackage{caption}
\usepackage{subcaption}
\usepackage{float}
\usepackage{xcolor}
\usepackage{pdfpages}

\newtheorem*{theorem*}{Theorem}
\newtheorem{theorem}{Theorem}[section]

\newtheorem{assumption}[theorem]{Assumption}
\newtheorem{corollary}[theorem]{Corollary}

\newtheorem{example}[theorem]{Example}
\newtheorem{lemma}[theorem]{Lemma}
\newtheorem*{lemma*}{Lemma}
\newtheorem*{definition*}{Definition}

\newtheorem{remark}[theorem]{Remark}

\numberwithin{equation}{section}
\numberwithin{figure}{section}
\usepackage[a4paper]{geometry}
\newcommand{\beq}{\begin{equation}}
\newcommand{\eeq}{\end{equation}}

\newcommand{\cA}{{\mathcal A}}

\newcommand{\cH}{{\cal H}}

\newcommand{\cD}{{\mathcal D}}

\newcommand{\bx}{\mathbf{x}}

\newcommand{\br}{\mathbf{r}}
\newcommand{\bV}{\mathbf{V}}

\newcommand{\bU}{\mathbf{U}}
\newcommand{\bW}{\mathbf{W}}

\newcommand{\bF}{\mathbf{F}}

\newcommand{\Com}{\mathbb{C}}

\newcommand{\ri}{{\rm i}}

\usepackage{color}
\usepackage{hyperref}
\definecolor{myblue}{rgb}{0,0,0.6}
\definecolor{darkgreen}{rgb}{0,0.5,0}
\hypersetup{colorlinks=true,
linkcolor=myblue,citecolor=myblue,filecolor=myblue,urlcolor=myblue}

\definecolor{escol}{rgb}{0,0,0.6}
\definecolor{sgcol}{rgb}{0,0,0.7}
\definecolor{estcol}{rgb}{0.5,0,0}
\definecolor{esnewcol}{rgb}{0,0.5,0}
\definecolor{lightgrayl}{RGB}{198,198,198}
\newcommand{\red}[1]{{\color{black}{#1}}}

\newcommand{\beqs}{\begin{equation*}}
\newcommand{\eeqs}{\end{equation*}}
\newcommand{\bit}{\begin{itemize}}
\newcommand{\eit}{\end{itemize}}
\newcommand{\ben}{\begin{enumerate}}
\newcommand{\een}{\end{enumerate}}
\newcommand{\bal}{\begin{align}}
\newcommand{\eal}{\end{align}}
\newcommand{\bals}{\begin{align*}}
\newcommand{\eals}{\end{align*}}
\newcommand{\bre}{\begin{remark}}
\newcommand{\ere}{\end{remark}}
\newcommand{\bpf}{\begin{proof}}
\newcommand{\epf}{\end{proof}}
\newcommand{\ble}{\begin{lemma}}
\newcommand{\ele}{\end{lemma}}
\newcommand{\bco}{\begin{corollary}}
\newcommand{\eco}{\end{corollary}}
\newcommand{\bex}{\begin{example}}
\newcommand{\eex}{\end{example}}
\newcommand{\bth}{\begin{theorem}}
\newcommand{\enth}{\end{theorem}}

\newcommand{\tfa}{\text{ for all }}
\newcommand{\tfor}{\text{ for }}

\newcommand{\ton}{\text{ on }}

\newcommand{\tand}{\text{ and }}

\newcommand{\Cdiso}{C_{{\rm dis}, N, 2}}
\newcommand{\Cdisoo}{C_{{\rm dis}, N, 1}}

\newcommand*{\N}[1]{\left\|#1\right\|}

\newcommand{\vertiii}[1]{{\left\vert\kern-0.25ex\left\vert\kern-0.25ex\left\vert #1 
    \right\vert\kern-0.25ex\right\vert\kern-0.25ex\right\vert}}

\newcommand{\matrixC}{\mathsf{C}}
\newcommand{\matrixI}{{\mathsf{I}}}

\usepackage[outdir=./]{epstopdf}
\usepackage{verbatim}

\usepackage{enumitem}

\newcommand{\matrixD}{\mathsf{D}}
\newcommand{\matrixA}{\mathsf{A}}

\newcommand{\matrixS}{\mathsf{S}}

\newcommand{\matrixM}{\mathsf{ M}}
\newcommand{\matrixB}{\mathsf{ B}}

\title{
Preconditioning FEM discretisations of the high-frequency \red{Helmholtz and} Maxwell equations by either perturbing the coefficients or adding absorption
}

\author{
E.~A.~Spence\thanks{Department of Mathematical Sciences, University of Bath, UK, {\tt e.a.spence@bath.ac.uk}}
}

\date{\today}

\begin{document}
\maketitle

\begin{abstract}
\red{This paper investigates the following question:~given a Galerkin matrix corresponding to a finite-element discretisation of either the Helmholtz or time-harmonic Maxwell equations with variable coefficients, suppose that the coefficients of the underlying PDE are perturbed; how good an approximate inverse (i.e., preconditioner) is the resulting Galerkin matrix to the original Galerkin matrix?}
An important special case is when the perturbation consists of adding absorption (in the spirit of ``shifted Laplacian preconditioning"). The results of this paper 
\red{improve the Helmholtz results in \cite{GaGrSp:15, GrPeSp:21} and extend these results to the time-harmonic Maxwell equations, confirming} 
a conjecture in the recent preprint \cite{LiHu:25}. 
\end{abstract}

\section{\red{Introduction}}

\subsection{Context and motivation}\label{sec:context}

The design of good preconditioners for the linear systems arising from finite-element discretisations of high-frequency time-harmonic wave problems -- such as the Helmholtz or time-harmonic Maxwell equations -- is a longstanding open problem; see, e.g., the review articles \cite{Er:08, ErGa:12, GaZh:19, GaZh:22} and the references therein.
This paper considers the theory of the following two, related, questions involving preconditioning these linear systems.

\paragraph{Preconditioning using absorption.}

The idea of preconditioning discretisations of the Helmholtz equation with (approximations of) discretisations of the same Helmholtz problem with added absorption was introduced in 
\cite{ErVuOo:04, ErOoVu:06}.
Seeking to rigorously understand this method, the paper \cite{GaGrSp:15} proved bounds on $\matrixI - \matrixA_\varepsilon^{-1}\matrixA$, where $\matrixA$ is the Helmholtz Galerkin matrix, and $\matrixA_\varepsilon$ is the Helmholtz Galerkin matrix with $k^2 \mapsto k^2 + \ri \varepsilon$, where $k$ is the (real and large) wavenumber. \cite{GaGrSp:15} proved that, for a particular nontrapping Helmholtz problem, $\|\matrixI - \matrixA_\varepsilon^{-1}\matrixA\|_2 \leq C \varepsilon/k$.
An important question is then whether the action of $\matrixA_\varepsilon^{-1}$ can be (provably) efficiently approximated when $\varepsilon\sim k$.
This question was investigated for one-level domain-decomposition methods in \cite{GSZ1}, and the recent work \cite{HuLi:24} proved that certain hybrid two-level Schwarz preconditioners (with problem-adapted basis functions in the coarse space) 
formed with absorption $\varepsilon\sim k$ 
give a  $k$-independent number of GMRES iterations when applied to the original Helmholtz problem. 
\paragraph{``Nearby preconditioning".}

The results of \cite{GaGrSp:15} were generalised in \cite{GrPeSp:21} to bounds on $\matrixI - \matrixA_2^{-1}\matrixA_1$, where $\matrixA_\ell$, $\ell=1,2$, are the Galerkin matrices of $k^{-2} \nabla \cdot(A_\ell\nabla u_\ell) + n_\ell u_\ell=f$. The motivation for studying this situation comes from uncertainty quantification (UQ):~calculating quantities
of interest for the solution of the Helmholtz equation 
$k^{-2} \nabla \cdot(A\nabla u) + n u=f$
with random coefficients requires the solution of many deterministic Helmholtz problems, each
with different coefficients $A$ and $n$. Bounds on $\matrixI - \matrixA_2^{-1}\matrixA_1$ then indicate to what extent
a previously-calculated inverse of one of the Galerkin matrices can be used as a preconditioner for other Galerkin matrices; see \cite[\S4.6]{Pe:20} and the recent work \cite{VaSc:25} for UQ algorithms using this idea.

\paragraph{Content and motivation for the present paper.}

The present paper 
\bit
\item proves the Helmholtz results of \cite{GaGrSp:15, GrPeSp:21} under less-restrictive assumptions than in \cite{GaGrSp:15, GrPeSp:21}, and 
\item generalises these results to the time-harmonic Maxwell equations. 
\eit 
Our main motivation is the recent work \cite{LiHu:25}:~this paper gives the first frequency-explicit analysis of a two-level domain-decomposition method for the high-frequency time-harmonic Maxwell equations, \emph{under the assumption that the Helmholtz results of \cite{GaGrSp:15, GrPeSp:21} hold also for the time-harmonic Maxwell equations}; the present paper shows that this assumption indeed holds.

\subsection{\red{Informal statement of the main result}}\label{sec:informal}

\red{
\paragraph{Set-up at the continuous level.}
Let $\cH$ be a Hilbert space and let $a_\ell(\cdot,\cdot):\cH\times \cH\to \Com$, $\ell=1,2,$ be sesquilinear forms that correspond to either the Helmholtz operator 
\beq\label{eq:Helmholtz}
-k^{-2} \nabla\cdot (\mu_\ell^{-1}\nabla )  - \epsilon_\ell
\eeq
or the time-harmonic Maxwell operator
\beq\label{eq:Maxwell}
k^{-2}{\rm curl}\big(\mu^{-1}_\ell {\rm curl}\big) - \epsilon_\ell
\eeq
with appropriate boundary conditions (encoded via the choices of $\cH$ and the sesquilinear form), and where the coefficients $\mu^{-1}_\ell$ and $\epsilon_\ell$ satisfy the natural assumptions to make the problems well-posed. 

Specific Helmholtz and Maxwell problems to which the main result applies are described in Lemmas \ref{lem:Helmholtz} and \ref{lem:Maxwell} below; we highlight that these include Helmholtz and Maxwell problems where the radiation condition is approximated by \emph{either} a perfectly-matched layer (PML) \emph{or} an impedance boundary condition.


\paragraph{Set-up at the discrete level.}
Given a finite-dimensional subspace $\cH_N\subset \cH$ (with $N= {\rm dim}(\cH_N)$), let $\matrixA_\ell, \ell=1,2,$ be the Galerkin matrices associated to $a_\ell(\cdot,\cdot)$. 
Let the matrix $\matrixD$ be such that 
\beq\label{eq:matrixD}
\N{v_N}_{\cH}^2 = \N{\bV}_{\matrixD}^2 \quad\tfor v_N=\sum_{j=1}^N V_j \phi_j\, \text{ with } \,V_j \in \mathbb{C}^N;
\eeq
i.e., $\|\cdot\|_{\matrixD}$ is the norm on $\Com^N$ inherited from the norm $\|\cdot\|_{\cH}$ on $\cH$. 
Let also $\|\cdot\|_{\matrixD}$ denote the induced norm on matrices. 

\paragraph{Informal statement of the main result (Theorem \ref{thm:main} below).}
The main result gives sufficient conditions under which $ \matrixI - \matrixA_2^{-1}\matrixA_1$ (measured in either $\|\cdot\|_{\matrixD}$ or the Euclidean norm $\|\cdot\|_2$) is controlled by $\|\mu_1-\mu_2\|_{L^\infty}$ and $\|\epsilon_1-\epsilon_2\|_{L^\infty}$ multiplied by the continuous inf-sup constant of $a_1(\cdot,\cdot)$.

In more detail, suppose that $\cH_N$ is sufficiently large so that the discrete inf-sup constant of $a_1(\cdot,\cdot)$ is bounded by a constant (independent of $N$) multiplied by the continuous inf-sup constant of $a_1(\cdot,\cdot)$; 
denote the latter by $\| \cA^{-1}_1\|_{\cH^*\to \cH}$ (see \S\ref{sec:recap} below for an explanation of this notation).
If 
\beqs
 \Big( \N{\mu^{-1}_1-\mu^{-1}_2}_{L^\infty} 
+ \N{\epsilon_1-\epsilon_2}_{L^\infty}\Big) C_1\big\| \cA^{-1}_1\big\|_{\cH^*\to \cH}\leq 1/2
\eeqs
then $ \matrixA_2^{-1}$ exists and 
\begin{align}\nonumber
&\max\Big\{
\big\| \matrixI - \matrixA_2^{-1}\matrixA_1 \big\|_{\matrixD}
, \big\| \matrixI - \matrixA_1 \matrixA_2^{-1}\big\|_{\matrixD^{-1}}\Big\}\\
&\hspace{3cm}\leq  2 \Big( \N{\mu^{-1}_1-\mu^{-1}_2}_{L^\infty} 
+ \N{\epsilon_1-\epsilon_2}_{L^\infty}\Big) C_1\big\| \cA^{-1}_1\big\|_{\cH^*\to \cH}.
\label{eq:main1a}
\end{align}
Furthermore, if $\mu_1=\mu_2$, then there exists $C_2$ (independent of $N$ and $k$) such that  
\beq\label{eq:main2a}
\max\Big\{
\big\| \matrixI - \matrixA_2^{-1}\matrixA_1 \big\|_{2},
\big\| \matrixI - \matrixA_1\matrixA_2^{-1} \big\|_2\Big\} \leq  C_2\N{\epsilon_1-\epsilon_2}_{L^\infty} C_1 \big\| \cA^{-1}_1\big\|_{\cH^*\to \cH},
\eeq
where $\|\cdot\|_2$ denotes the Euclidean norm on matrices (induced by the Euclidean norm on vectors).
}

\red{
\subsection{Discussion}

\subsubsection{The main result specialised to preconditioning with absorption}
In the set-up of \S\ref{sec:informal}, let $\mu_2=\mu_1$ and $\epsilon_2= (1+\ri \alpha) \epsilon_1$ with $\alpha>0$ (i.e., $a_2(\cdot,\cdot)$ is formed by addition absorption to $a_1(\cdot,\cdot)$). 
If 
\beqs
\alpha C_1\big\| \cA^{-1}_1\big\|_{\cH^*\to \cH}\leq 1/2
\eeqs
(i.e., the absorption is sufficiently small, depending on the inf-sup constant of $a_1(\cdot,\cdot)$),
then $ \matrixA_2^{-1}$ exists,
\begin{align}
&\max\Big\{
\big\| \matrixI - \matrixA_2^{-1}\matrixA_1 \big\|_{\matrixD}
, \big\| \matrixI - \matrixA_1 \matrixA_2^{-1}\big\|_{\matrixD^{-1}}\Big\} \leq 2 \alpha  C_1\big\| \cA^{-1}_1\big\|_{\cH^*\to \cH}
\label{eq:abs1a}
\end{align}
and
\beq\label{eq:abs2a}
\max\Big\{
\big\| \matrixI - \matrixA_2^{-1}\matrixA_1 \big\|_{2},
\big\| \matrixI - \matrixA_1\matrixA_2^{-1} \big\|_2\Big\} \leq  C_2\alpha C_1 \big\| \cA^{-1}_1\big\|_{\cH^*\to \cH}.
\eeq

\subsubsection{Recovering the Helmholtz results of \cite{GaGrSp:15} about preconditioning with absorption}
\cite{GaGrSp:15} considers preconditioning the Helmholtz equation with absorption by letting $k^2\mapsto k^2 + \ri \varepsilon$. In the set-up above this corresponds to setting $\alpha:=\varepsilon/k^2$. By Remark \ref{rem:Csol} below, $\big\| \cA^{-1}_1\big\|_{\cH^*\to \cH}\sim k$ when the Helmholtz problem is nontrapping. The right-hand side of \eqref{eq:abs2a} is then proportional to $\varepsilon/k$, which is the result proved in  \cite[Theorem 1.4]{GaGrSp:15} for the specific case of the Helmholtz interior impedance problem.

\subsubsection{The condition that the discrete inf-sup constant is bounded by the continuous inf-sup constant}\label{sec:final_day}

As stated in \S\ref{sec:informal}, the main result holds under the assumption that the dimension of the finite-dimensional approximation space is large enough so that the discrete inf-sup constant is bounded by a constant multiple of the continuous inf-sup constant.

This condition on the discrete inf-sup constant has been proven to hold for Helmholtz problems in the so-called \emph{asymptotic regime}, i.e., when the finite-dimensional space is sufficiently large so that the Galerkin solution is quasi-optimal; see \cite[Theorem 4.2]{MeSa:10} and the discussion in Remark \ref{rem:link} below. Recall that, for the $h$-version FEM applied to Helmholtz non-trapping problems, the asymptotic regime is when $(kh)^p k$ is sufficiently small \cite{MeSa:10, ChNi:20}.

In \S\ref{sec:is} we prove a new result showing that the discrete inf-sup constant is bounded by a constant multiple of the continuous inf-sup constant in the \emph{preasymptotic regime} when, for the $h$-version FEM applied to Helmholtz non-trapping problems, $(kh)^{2p}k$ is sufficiently small \cite{GS3}. 

Given a Helmholtz or Maxwell problem of interest, the ``recipe" to use the main result of this paper is therefore the following:~
\ben
\item Show that the problem satisfies Assumption \ref{ass:abstract} below. Note that Lemmas \ref{lem:Helmholtz} and \ref{lem:Maxwell} show that many commonly-studied Helmholtz and Maxwell problems satisfy this assumption.
\item Determine sufficient conditions for the discrete inf-sup constant to be bounded by a constant multiple of the continuous inf-sup constant -- these conditions will involve regularity assumptions on the domain and coefficients, as well as conditions on the dimension of the approximation space. Section \ref{sec:is} below details these conditions (with references) for the $h$-FEM applied to the Helmholtz and Maxwell problems in Lemmas \ref{lem:Helmholtz} and \ref{lem:Maxwell}.
\een 


\subsubsection{How the results of this paper improve the Helmholtz results of \cite{GaGrSp:15, GrPeSp:19}}

Remark \ref{rem:improve} below discusses in detail how the results of this paper improve the Helmholtz results of \cite{GaGrSp:15, GrPeSp:19}. We highlight here that the main way (aside from extending the results to the time-harmonic Maxwell equations) is in proving that bounds \eqref{eq:main1a}-\eqref{eq:abs2a} hold if  the discrete inf-sup constant is bounded by a constant multiple of the continuous inf-sup constant -- the analyses in \cite{GaGrSp:15, GrPeSp:19} require more restrictive conditions on the finite-dimensional space.


\subsubsection{Numerical experiments investigating the  sharpness of the bounds \eqref{eq:main1a}, \eqref{eq:main2a}, \eqref{eq:abs1a}, and \eqref{eq:abs2a}}

\cite[\S5]{GaGrSp:15} and \cite[\S2.2.2]{GrPeSp:21} contain numerical experiments investigating the sharpness of the bounds \eqref{eq:abs2a} and \eqref{eq:main1a}/\eqref{eq:main2a} respectively.
We highlight that these experiments consider piecewise-linear FEM discretisations in the preasymptotic regime, i.e., $(kh)^{2p}k$ with $p=1$; these FEM discretisations were not covered by the theory in \cite{GaGrSp:15, GrPeSp:21} (see the discussion in \S\ref{sec:final_day} and in Remark \ref{rem:improve} below) but are covered by the results of the present paper.

\subsubsection{The use of the bound \eqref{eq:abs1a} in the DD analyses of \cite{HuLi:24} and \cite{LiHu:25}}
The Helmholtz two-level DD analysis of \cite{HuLi:24} assumes that the bound \eqref{eq:abs1a} holds for the Helmholtz interior impedance problem 
(see the assumptions of \cite[Corollary 6.1]{HuLi:24}).
Similarly, the Maxwell two-level DD analysis of \cite{LiHu:25} assumes that the bound \eqref{eq:abs1a} holds for the Maxwell interior impedance problem (see the discussion on \cite[Page 2]{LiHu:25}). 
The main result of this paper, Theorem \ref{thm:main}, coupled with results about the discrete inf-sup constant in \S\ref{sec:is} below then give sufficient conditions on the finite-element spaces for these assumptions to be satisfied.

}

\paragraph{Outline of the \red{rest of the} paper.}

\S\ref{sec:assumptions} states the abstract assumptions under which the main results are proved 
and describes Helmholtz and Maxwell problems that are covered by these abstract assumptions.
\S\ref{sec:recap} recaps properties of the inf-sup constant of a sesquilinear form and its $k$-dependence for Helmholtz and Maxwell.
\S\ref{sec:main} states the main results, and \S\ref{sec:proofs} proves them.

\section{
\red{The class of Helmholtz and Maxwell problems considered}
}\label{sec:assumptions}

\red{The main result is proved under the following abstract assumptions.}

\begin{assumption}[\red{Abstract assumptions}]\label{ass:abstract}
$\cH \subset \cH_0$ are Hilbert spaces with $\|\cdot\|_{\cH_0}\leq \|\cdot\|_{\cH}$ \red{and}
$\cH_0$ \red{identified} with its dual so that $\cH\subset \cH_0\subset\cH^*$. 
The linear operator $\cD:\cH\to \cH_0$ with $\|\cD\|_{\cH\to \cH_0}\leq 1$.
$b(\cdot,\cdot)$ is a continuous sesquilinear form on $\cH$, 
\beq\label{eq:sesqui}
a_\ell(u,v):=  \big( \mu_\ell^{-1} \cD u, \cD v\big)_{\cH_0} + b(u,v) - \big(\epsilon_\ell u,v\big)_{\cH_0}, \quad\ell=1,2,
\eeq
where $\mu_\ell^{-1}:\cH_0\to \cH_0$ and $\epsilon_\ell:\cH_0\to \cH_0$ are bounded linear operators for $\ell=1,2$. 
Finally for $\ell=1,2,$ there exists $C_{\rm G1,\ell}, C_{\rm G2,\ell}>0$ such that 
\beq\label{eq:Ccoer}
\big| a_\ell(v,v) +C_{\rm G2,\ell}\N{v}^2_{\cH_0} \big|  \geq C_{\rm G1,\ell} \N{v}^2_{\cH}\quad\tfa v\in \cH.
\eeq
\end{assumption}

We make three \red{immediate} remarks.
\bit
\item 
These assumptions imply that $a_\ell(\cdot,\cdot)$ is continuous; i.e., there exist $C_{\rm cont,\ell}>0$ such that
\beqs
\big| a_\ell(u,v)\big| \leq C_{\rm cont,\ell} \N{u}_{\cH} \N{v}_{\cH} \quad\tfa u,v\in \cH.
\eeqs
\item If $a_{\red \ell}(\cdot,\cdot)$ satisfies the G\aa rding inequality 
\beq\label{eq:Garding}
\Re a_\ell(v,v)\geq C_{\rm G1,\ell} \N{v}^2_{\cH} - C_{\rm G2,\ell} \N{v}^2_{\cH_0} \quad\tfa v\in \cH,
\eeq
then \eqref{eq:Ccoer} holds (since $|z|\geq \Re z$ for all $z\in\Com$).
\item 
 If $a_\ell(\cdot,\cdot)$ satisfies these assumptions, then so does 
 \beqs
 a_\ell^*(u,v):= \overline{a_\ell(v,u)}= \big( \mu_\ell^* \cD u, \cD v\big)_{\cH_0} + \overline{b(v,u)} - \big(\epsilon_\ell^* u,v\big)_{\cH_0},
 \eeqs
  with the same $C_{\rm G1, \ell}, C_{\rm G2, \ell}$ (this is because 
  $| a_\ell(v,v) +C_{\rm G2,\ell}\N{v}^2_{\cH_0}|= \big|\overline{a_\ell(v,v)} +C_{\rm G2,\ell}\N{v}^2_{\cH_0} \big|$).
\eit

\red{We now give examples of Helmholtz and Maxwell problems satisfying Assumption \ref{ass:abstract}.}

\red{
\begin{lemma}[Helmholtz problems satisfying Assumption \ref{ass:abstract}]\label{lem:Helmholtz}
Let $\Omega$ be a bounded Lipschitz open set, $\cH_0 = L^2(\Omega)$ (scalar valued), $\cD:= k^{-1} \nabla$, and
\beq\label{eq:norm}
\N{v}^2_{\cH}:= \N{k^{-1}\nabla v}^2_{L^2(\Omega)} + \N{v}^2_{L^2(\Omega)}.
\eeq
For $\ell=1,2,$ let $\mu_\ell^{-1}$ be bounded, symmetric, matrix functions on $\Omega$ with ${\rm essinf}_\Omega\Re \mu_\ell^{-1}>0$ (in the sense of quadratic forms) and $\epsilon_\ell$ be bounded, scalar-valued functions on $\Omega$. 

If one of the following three points holds, then Assumption \ref{ass:abstract} holds.

(i) 
\beqs
\cH:= \big\{ u\in H^1(\Omega) \,:\, \text{ $u=0$ on at least one connected component of $\partial\Omega$}\big\}
\eeqs
and $b(\cdot,\cdot)=0$. 

(ii)  With $\Gamma_{\rm imp}$  a non-empty connected component of $\partial \Omega$, and $\Gamma_{\rm D}$ a (possibly empty) connected component of $\partial \Omega$ that is not equal to $\Gamma_{\rm imp}$,
\beqs
\cH:= \big\{ v\in H^1(\Omega)\,:\, v=0 \ton \Gamma_{\rm D}\big\},
\eeqs
and
\beqs
b(u,v):= -\ri k^{-1} \big( \theta u, v\big)_{L^2(\Gamma_{\rm imp})}
\eeqs
for some real-valued $\theta\in L^\infty(\Gamma_{\rm imp})$ such that ${\rm essinf}_{\Gamma_{\rm imp}}\theta >0$. 

(iii) 
$\Omega= \{ x: |x|\leq R\}\setminus \overline{\Omega_-}$, where $\Omega_-$ is a bounded Lipschitz open set with connected open complement (so that the scattering problem with obstacle $\Omega_-$ makes sense). 
$\cH$ equals either $H^1(\Omega)$ or $\{ v\in H^1(\Omega) \,:\, v=0 \ton \partial\Omega_-\}$.
\beqs
b(u,v):= -\ri k^{-1} \big( {\rm DtN}_k u, v\big)_{L^2(\Gamma_{\rm imp})},
\eeqs
where ${\rm DtN}_k$ is the exact Dirichlet-to-Neumann map for the Helmholtz equation in the exterior of $ \{ x: |x|\leq R\}$; see, e.g., \cite[Equations 3.5 and 3.6]{ChMo:08} for explicit expressions in terms of Fourier series (in 2-d) or spherical harmonics (in 3-d) and Bessel and Hankel functions. 
\end{lemma}
}

\red{
\bre[The Helmholtz problems covered in Lemma \ref{lem:Helmholtz}]
The fact that  $\cD:= k^{-1} \nabla$ implies that the sesquilinear form \eqref{eq:sesqui} corresponds 
 to the Helmholtz operator \eqref{eq:Helmholtz}
Parts (i), (ii), and (iii) of Lemma \ref{lem:Helmholtz} then cover three common ways of approximating the Sommerfeld radiation condition (i.e., the fact the scattered wave travels ``outward" from the scatterer; see, e.g., \cite[\S3.2]{CoKr:83}). 

Indeed, in Part (i), zero Dirichlet boundary conditions are imposed on at least one connected component of $\partial\Omega$, with then zero Neumann boundary conditions imposed on the rest.
The requirements on $\mu_\ell$ and $\epsilon_\ell$ in Lemma \ref{lem:Helmholtz} are then satisfied by the coefficients coming from a radial perfectly-matched layer (PML) by, e.g., \cite[Lemma 2.3]{GLSW1}.

In Part (ii), the sesquilinear form corresponds to the Helmholtz operator \eqref{eq:Helmholtz} with the impedance boundary condition
\beqs
k^{-1}\partial_{n, \mu^{-1}} u - \ri \theta u =g_{\rm imp} \quad\ton \Gamma_{\rm imp}
\eeqs
$u=0$ on $\Gamma_{\rm D}$, and $k^{-1}\partial_{n, \mu^{-1}} u=g_{\rm N}$ on $\partial\Omega \setminus (\Gamma_{\rm D} \cup \Gamma_{\rm imp})$
(where the data $g_{\rm imp}$ and $g_{\rm N}$ depend on the given right-hand side in the variational problem). 
Here $\partial_{n, \mu^{-1}}$ is the conormal derivative; recall that this is such that $\partial_{n, \mu^{-1}} u := n \cdot \mu^{-1}\nabla u$ when $u\in H^2(\Omega)$, where $n$ is the outward-pointing unit normal vector to $\partial\Omega$ (see, e.g., \cite[Lemmas 4.2 and 4.3]{Mc:00}).

In Part (iii), the sesquilinear form corresponds to the Helmholtz operator \eqref{eq:Helmholtz} with the exact Dirichlet-to-Neumann map imposed on the boundary of a ball.
\ere
}

\red{
\bpf[Proof of Lemma \ref{lem:Helmholtz}]
We need to check that (a) $\|\cD\|_{\cH\to \cH_0}\leq 1$, (b) $b(\cdot,\cdot)$ is continuous, and (c) the inequality \eqref{eq:Ccoer} holds.

Regarding (a):~this follows immediately from the definition of the norm \eqref{eq:norm} and the fact that 
 $\cD:= k^{-1} \nabla$.
 
 Regarding (b):~for Part (i), $b(\cdot,\cdot)=0$ and so is automatically continuous. For Part (ii), continuity of $b(\cdot,\cdot)$ follows from the standard 
 multiplicative trace inequality (see, e.g., \cite[
Theorem 1.5.1.10, last formula on p. 41]{Gr:85}, \cite[Theorem 1.6.6]{BrSc:08}). For Part (iii), continuity of $b(\cdot,\cdot)$ follows from, e.g., \cite[Lemma 3.3, Part 1]{MeSa:10}. 
 
Regarding (c):~for Part (i), the inequality \eqref{eq:Ccoer} follows immediately since $b(\cdot,\cdot)=0$. For Part (ii), $a(\cdot,\cdot)$ satisfies the G\aa rding inequality \eqref{eq:Garding} (with 
the term on $\Gamma_{\rm imp}$ playing no role because of the real part on the left-hand side of \eqref{eq:Garding}), and thus \eqref{eq:Ccoer} holds. For Part (iii), recall that $\Re (-\ri {\rm DtN}_k)\geq 0$ by, e.g., \cite[Corollary 3.1]{ChMo:08} or \cite[Lemma 3.3, Part 2]{MeSa:10}, so that $a(\cdot,\cdot)$ satisfies the G\aa rding inequality \eqref{eq:Garding}.
\epf 
}

\red{
\begin{lemma}[Maxwell problems satisfying Assumption \ref{ass:abstract}]\label{lem:Maxwell}
Let $\Omega$ be a bounded Lipschitz open set, $\cH_0=L^2(\Omega)$ (vector valued), and 
$\cD :=k^{-1}{\rm curl}$.

For $\ell=1,2,$ let $\mu_\ell^{-1}$ and $\epsilon_\ell$ be bounded, symmetric, matrix functions on $\Omega$ with $\red{{\rm essinf}_\Omega}\Re \mu_\ell^{-1}>0$ and $\red{{\rm essinf}_\Omega}\Re \epsilon_\ell>0$ (in the sense of quadratic forms). 

If one of the following two points holds, then Assumption \ref{ass:abstract} holds.

(i) 
\beqs
\cH := H_0({\rm curl};\Omega) = \big\{ v \in H({\rm curl};\Omega) : v_T =0\big\}
\eeqs 
where $v_T$ is the tangential trace such that $v_T:= (n\times v)\times n$ for smooth $v$, where $n$ is the outward-pointing unit normal vector (see, e.g., \cite[Equation 3.46 and Theorem 3.31]{Mo:03}).
\beq\label{eq:normM1}
\N{v}^2_{\cH}:= \N{k^{-1}{\rm curl}\, v}^2_{L^2(\Omega)} + \N{v}^2_{L^2(\Omega)}
\eeq
and $b(\cdot,\cdot)=0$. 

(ii) With $\Gamma_{\rm imp}$  a non-empty connected component of $\partial \Omega$, 
\beqs
\cH:= \big\{ u\in H^1(\Omega) \,:\, \text{ $u_T=0$ on $\partial\Omega\setminus \Gamma_{\rm imp}$ and $u_T \in L^2(\Gamma_{\rm imp})$}\big\},
\eeqs
\beq\label{eq:normM2}
\N{v}^2_{\cH}:= k^{-2} \N{{\rm curl}\, v}^2_{L^2(\Omega)} + \N{v}^2_{L^2(\Omega)}+ k^{-1}\N{v_T}^2_{L^2(\Gamma_{\rm imp})},
\eeq
and
\beqs
b(u,v):= -\ri k^{-1} \big( \theta u_T, v_T\big)_{L^2(\Gamma_{\rm imp})},
\eeqs
for some real-valued $\theta\in L^\infty(\Gamma_{\rm imp})$ such that ${\rm essinf}_{\Gamma_{\rm imp}}\theta >0$. In addition, $\mu_\ell^{-1}$ and $\epsilon_\ell$ are real-valued.
\end{lemma} 
}

\red{
\bre[The Maxwell problems covered in Lemma \ref{lem:Helmholtz}]
The fact that  $\cD:= k^{-1} {\rm curl}$ implies that the sesquilinear form \eqref{eq:sesqui} corresponds 
 to the Maxwell operator 
and the sesquilinear form corresponds to the Maxwell operator \eqref{eq:Maxwell}
Parts (i) and (ii) of Lemma \ref{lem:Helmholtz} then cover two common ways of approximating the  radiation condition.

Indeed, the requirements on $\mu_\ell$ and $\epsilon_\ell$ in Part (i) are satisfied by the coefficients coming from a radial PML by, e.g., \cite[Appendix A]{CGS1}.
In Part (ii), the sesquilinear form corresponds to the Maxwell operator \eqref{eq:Maxwell} with the impedance boundary condition
\beqs
k^{-1}\big(\mu_\ell^{-1} {\rm curl}\, u\big)\times n - \ri \theta u_T =g \quad\ton \Gamma_{\rm imp}
\eeqs
(where $g$ depends on the given right-hand side) and the PEC boundary condition $u_T=0$ on $\Gamma\setminus \Gamma_{\rm imp}$.
\ere
}

\red{
\bpf[Proof of Lemma \ref{lem:Maxwell}]
As in the Helmholtz case, we need to check that (a) $\|\cD\|_{\cH\to \cH_0}\leq 1$, (b) $b(\cdot,\cdot)$ is continuous, and (c) the inequality \eqref{eq:Ccoer} holds.

Regarding (a):~this follows immediately from the definitions of the norm \eqref{eq:normM1}/\eqref{eq:normM2} and the fact that 
 $\cD:= k^{-1} {\rm curl}$.

Regarding (b):~for Part (i), continuity of $b(\cdot,\cdot)$ is immediate (since it is zero). For Part (ii), continuity of $b(\cdot,\cdot)$ follows from the Cauchy--Schwarz inequality and the definition of the norm \eqref{eq:normM2}.

Regarding (c):~for Part (i), the inequality \eqref{eq:Ccoer} follows immediately since $b(\cdot,\cdot)=0$. For Part (ii), since $\mu^{-1}_\ell$ and $\epsilon_\ell$ are both real (by assumption), $\Im (a(v,v))$ consists only of $\Im (b(v,v))$, which is bounded below by $k^{-1} ({\rm essinf}_{\Gamma_{\rm imp}} \theta )\|v\|^2_{L^2(\Gamma_{\rm imp})}$; the bound \eqref{eq:Ccoer} then follows.
\epf
}

\bre[The reason for the particular weighting with $k$]
Most papers on the numerical analysis of the Helmholtz and time-harmonic Maxwell equations write these equations as $(\Delta +k^2) u=f$ and $({\rm curl}\, {\rm curl}  -k^2) E=f$ (in the constant coefficient case) and use the norms $\|v\|^2_{H^1} = \|\nabla v\|^2_{L^2} + k^2\|v\|^2_{L^2}$ and $\|w\|_{H({\rm curl})}^2 = \|{\rm curl} w\|^2_{L^2} + k^2 \|w\|^2_{L^2}$. 
The advantage of rescaling the PDEs and norms so that every derivative appears with a $k^{-1}$ is that none of the constants in \S\ref{sec:assumptions} (i.e. $C_{\rm G1, \ell}, C_{\rm G2, \ell},$ and $C_{\rm cont, \ell}$) depend on $k$, and (consequentially) the $k$-dependence of the norm of the solution operator is the same between any two spaces for which the solution operator is defined; i.e., the $\cH_0\to \cH_0$ norm, the $\cH_0\to \cH$ norm, and the $\cH^*\to \cH$ norm all have the same $k$-dependence -- see Lemma \ref{lem:is1} below.
\ere

\section{Recap of properties of the inf-sup constant of $a_\ell(\cdot,\cdot)$}\label{sec:recap}

\red{
In this section, we recap properties of the inf-sup constant of $a_\ell(\cdot,\cdot)$, since it appears in the statement of the main result (Theorem \ref{thm:main} below).}

Let $\cA_\ell:\cH\to \cH^*$ be the operator associated to $a_\ell(\cdot,\cdot)$; i.e., $\langle \cA_\ell u,v \rangle_{\cH^*\times \cH}= a_\ell(u,v)$ for all $u,v\in \cH$. 
We recall the following standard result.

\begin{theorem}[Inf-sup condition equivalent to bound on inverse operator]\label{thm:inverse}
The conditions 
\beq\label{eq:is_cond}
\inf_{u\in \cH\setminus\{0\}} \sup_{v\in \cH\setminus\{0\}} \frac{
\big| a_\ell(u,v)\big|
}{
\N{u}_{\cH} \N{v}_{\cH}
}
\geq 
\gamma^{-1}
\text{ and for all } v\in \cH\setminus \{0\}, \sup_{u\in \cH\setminus\{0\}} |a_\ell(u,v)|>0
\eeq
and
\beqs
\big\|\cA_\ell^{-1}\big\|_{\cH^*\to \cH} \leq \gamma
\eeqs
are equivalent.
\end{theorem}

\bpf[References for the proof]
See, e.g., \cite[Lemma 6.5.3]{Ha:92}, \cite[Theorem 2.1.44]{SaSc:11}. The first bound in \eqref{eq:is_cond} is that $\|\cA_\ell u\|_{\cH^*}\geq \gamma^{-1}\|u\|_{\cH}$ for all $u\in \cH$, and the second bound in \eqref{eq:is_cond} implies that the kernel of $\cA_\ell^*$ is empty, and thus the image of $\cA_\ell$ equals $\cH^*$; see, e.g., \cite[Theorem 2.13(iii)]{Mc:00}.
\epf

The following lemma, proved in \S\ref{sec:Geneva} below using the G\aa rding-type inequality \eqref{eq:Ccoer}, 
shows that all norms of $\cA_\ell^{-1}$ are equivalent, with constants depending only on $C_{\rm G1,\ell}$ and $C_{\rm G2,\ell}$ (and hence independent of $k$).

\begin{lemma}[$k$-independent equivalence of norms of $\cA_\ell^{-1}$] \label{lem:is1}
If $\cA_\ell^{-1}:\cH^*\to \cH$ exists, then
\beq\label{eq:sunny1}
\big\|\mathcal{A}_\ell^{-1}\big\|_{\cH_0\to \cH} \leq \big\|\mathcal{A}_\ell^{-1}\big\|_{\cH^*\to \cH} \leq (C_{\rm G1, \ell})^{-1} \Big( 1 + C_{\rm G2, \ell}\big\|\mathcal{A}_\ell^{-1}\big\|_{\cH_0 \to \cH}\Big)
\eeq
and 
\beq\label{eq:sunny2}
\big\|\mathcal{A}_\ell^{-1}\big\|_{\cH_0\to \cH_0} \leq \big\|\mathcal{A}_\ell^{-1}\big\|_{\cH_0\to \cH} \leq 
(C_{\rm G1,\ell})^{-1/2}\big\|\mathcal{A}_\ell^{-1}\big\|_{\cH_0\to \cH_0}\sqrt{ C_{\rm G2, \ell}  + \big\|\mathcal{A}_\ell^{-1}\big\|_{\cH_0\to \cH_0}^{-1}}.
\eeq
\end{lemma}

\bre[The $k$-dependence of $\|\cA^{-1}_\ell\|$ for Helmholtz and Maxwell]\label{rem:Csol}
When $\mu_\ell$ and $\epsilon_\ell$ are both constant multiples of the identity in part of the domain, $\|\cA^{-1}_\ell\| \geq Ck $ -- this can be proved 
by considering data that is a cut-off function multiplied by a plane wave; see, e.g., \cite[Lemma 3.10]{ChMo:08} (for Helmholtz) and \cite[\S1.4.1]{ChMoSp:23}, \cite[Example 3.4]{MeSa:23} (for Maxwell). 

For the scattering problem, if the scatterer is nontrapping then $\|\cA^{-1}_\ell\| \leq Ck $; this has been proved in a wide variety of Helmholtz cases (see, \red{e.g.,}, the literature review in \cite{GrPeSp:19}) and for the Maxwell problem with certain nontrapping coefficients in 
\cite{ChMoSp:23} or a nontrapping PEC obstacle in \cite[\S2]{Ya:88}. Under trapping, $\|\cA^{-1}_\ell\|\gg k$; see the results and literature reviews in \cite{LSW1}, \cite{ChSpGiSm:20}.

For the PML problem, \cite[Theorem 1.6]{GLS2} proved that the norm of the solution operator of the Helmholtz radial PML problem is bounded by the norm of the solution operator of the corresponding Helmholtz scattering problem; we expect that the same result holds for the Maxwell PML problem. Indeed,  this result was recently proved \red{-- up to a factor of $k$ loss --} for the case of no scatterer and Cartesian PML in \cite[Lemma 10]{CuWaXi:25}.
\ere

\section{Statement of the main results}\label{sec:main}

\red{This section states precisely the main results of the paper; the proofs of these results are then given in \S\ref{sec:proofs}.}

\subsection{Notation for the Galerkin matrices}

Let $\cH_N\subset \cH$ be a finite-dimensional space with basis $\{\phi_j\}_{j=1}^N$. Let 
\beq\label{eq:matrixA}
\matrixA_\ell = \matrixS_{\mu^{-1}_\ell} + \matrixB - \matrixM_{\epsilon_\ell}, \quad\ell=1,2,
\eeq
where
\beq\label{eq:matrixS}
(\matrixS_{\mu^{-1}_\ell})_{ij} = \big( \mu^{-1}_\ell\cD \phi_j ,\cD\phi_i\big)_{\cH_0}, \quad (\matrixM_{\epsilon_\ell})_{ij} = \big(\epsilon_\ell \phi_j ,\phi_i\big)_{\cH_0}, \quad\tand\quad \matrixB_{ij} = b(\phi_j,\phi_i).
\eeq
Let the matrix $\matrixD$ be such that \red{the norm relation} \eqref{eq:matrixD} holds.
We also use the notation $\|\cdot\|_{\matrixD}$ to denote the induced norm on matrices. Finally, let $\|\cdot\|_2$ denote the Euclidean norm on $\Com^N$, and let $m_{\pm}$ be such that 
\beq\label{eq:normequiv}
m_-\|\bV\|_2 \leq \| v_N\|_{\cH_0} \leq m_+ \|\bV\|_2\quad\tfa v_N\in \cH_N.
\eeq

\bre
With $\matrixM$ the mass matrix, i.e., $(\matrixM)_{ij}= ( \phi_j ,\phi_i)_{\cH_0}$, \eqref{eq:normequiv} is equivalent to the bounds
\beqs
m_-^2\|\bV\|_2^2 \leq \big( \matrixM \bV, \bV\big)_2 \leq m_+^2 \|\bV\|_2^2\quad\tfa \bV\in \Com^N;
\eeqs
i.e., the quantity $(m_+/m_-)^2$  (whose square root appears in the bound \eqref{eq:main2} below) is the ratio of the largest to the smallest eigenvalue -- i.e., the condition number -- of the mass matrix $\matrixM$.
\red{For the $h$-FEM with quasi-uniform meshes, the ratio $m_+/m_-$ is independent of $h$; see, e.g., \cite[Equation 4.2]{GaGrSp:15}, \cite[Lemma 5.1]{GrPeSp:21}, \cite[Lemma 4.6]{Pe:19}.}
\ere

\subsection{The main result:~bounds on $\| \matrixI - \matrixA_2^{-1}\matrixA_1\|$
and $\| \matrixI -\matrixA_1 \matrixA_2^{-1}\|$}

\begin{theorem}[Bounds on $\| \matrixI - \matrixA_2^{-1}\matrixA_1\|$ and $\| \matrixI -\matrixA_1 \matrixA_2^{-1}\|$]
\label{thm:main}
Suppose that the assumptions of \S\ref{sec:assumptions} hold. 
Suppose that $(\cH_N)_{N=1}^\infty$ are such that there exists $C_1>0$ such that, for all $N\in \mathbb{Z}^+$,
\begin{subequations}\label{eq:dis_both}
\beq\label{eq:dis}
\inf_{u_N\in \cH_N\setminus\{0\}} \sup_{v_N\in \cH_N\setminus\{0\}} \frac{
\big| a_1(u_N,v_N)\big|
}{
\N{u_N}_{\cH} \N{v_H}_{\cH}
}
\geq 
\frac{1}{C_1 \big\| \cA_1^{-1}\big\|_{\cH^*\to \cH}}
\eeq
and
\beq\label{eq:dist}
\tfa v_N\in \cH_N\setminus \{0\},\, \sup_{u_N\in \cH_N\setminus\{0\}} |a_1(u_N,v_N)|>0
\eeq
\end{subequations}
(i.e., the discrete inf-sup constant of $a_1(\cdot,\cdot)$ is bounded below by a constant times the continuous inf-sup constant \eqref{eq:is_cond}).
If 
\beq\label{eq:FridaySun1}
 \Big( \N{\mu^{-1}_1-\mu^{-1}_2}_{\cH_0\to \cH_0} 
+ \N{\epsilon_1-\epsilon_2}_{\cH_0\to \cH_0}\Big) C_1\big\| \cA^{-1}_1\big\|_{\cH^*\to \cH}\leq 1/2
\eeq
then $ \matrixA_2^{-1}$ exists and 
\begin{align}\nonumber
&\max\Big\{
\big\| \matrixI - \matrixA_2^{-1}\matrixA_1 \big\|_{\matrixD}
, \big\| \matrixI - \matrixA_1 \matrixA_2^{-1}\big\|_{\matrixD^{-1}}\Big\}\\
&\hspace{3cm}\leq  2 \Big( \N{\mu^{-1}_1-\mu^{-1}_2}_{\cH_0\to \cH_0} 
+ \N{\epsilon_1-\epsilon_2}_{\cH_0\to \cH_0}\Big) C_1\big\| \cA^{-1}_1\big\|_{\cH^*\to \cH}.
\label{eq:main1}
\end{align}
Furthermore, if $\mu_1=\mu_2$, then 
\beq\label{eq:main2}
\max\Big\{
\big\| \matrixI - \matrixA_2^{-1}\matrixA_1 \big\|_{2},
\big\| \matrixI - \matrixA_1\matrixA_2^{-1} \big\|_2\Big\} \leq  2\frac{m_+}{m_-}\N{\epsilon_1-\epsilon_2}_{\cH_0\to \cH_0} C_1 \big\| \cA^{-1}_1\big\|_{\cH^*\to \cH}.
\eeq
\end{theorem}

We make the following four immediate remarks.
\bit
\item[(i)] The assumptions of Theorem \ref{thm:main} -- by design -- involve  $a_1(\cdot,\cdot)$ and not $a_2(\cdot,\cdot)$; i.e., the assumptions are on the problem that we want to solve (involving $\matrixA_1$), rather than the problem used for preconditioning (involving $\matrixA_2$).
\item[(ii)] Since $ \matrixI - \matrixA_2^{-1}\matrixA_1=  \matrixA_2^{-1}( \matrixA_2-\matrixA_1)$, 
the right-hand sides of both \eqref{eq:main1} and \eqref{eq:main2} should be understood as $\|\matrixA_2-\matrixA_1\| \|\matrixA_2^{-1}\|$; i.e., the norm of the perturbation times the norm of the solution operator. 
The key point is that the conditions \eqref{eq:dis} and \eqref{eq:FridaySun1} mean that $ \|\matrixA_2^{-1}\|$ is bounded by $\|\cA_2^{-1}\|_{\cH^*\to \cH}$, which in turn is bounded by $\|\cA_1^{-1}\|_{\cH^*\to \cH}$.
\item[(iii)]
\S\ref{sec:is} 
 gives sufficient conditions for \eqref{eq:dis} to hold (including a new result on this); i.e., conditions under which the discrete inf-sup constant is bounded by \red{a constant times} the continuous inf-sup constant.
\item[(iv)] In all our Helmholtz and Maxwell examples, 
the $\cH_0\to \cH_0$ norms of the coefficient differences appearing on the right-hand sides of both \eqref{eq:main1} and \eqref{eq:main2} are bounded by their $L^\infty$ norms (since $\|\mu v\|_{L^2(\Omega)}\leq \|\mu\|_{L^\infty(\Omega)}\N{v}_{L^2(\Omega)}$ for all $v\in L^2(\Omega)$).
\eit

\begin{corollary}[Theorem \ref{thm:main} specialised to preconditioning with absorption]
\label{cor:absorption}
Suppose that the assumptions of \S\ref{sec:assumptions} hold with $\mu_2=\mu_1$ and $\epsilon_2= (1+\ri \alpha) \epsilon_1$ with $\alpha>0$. 
Suppose that $(\cH_N)_{N=1}^\infty$ are such that there exists $C_1>0$ such that \eqref{eq:dis_both} hold for all $N\in \mathbb{Z}^+$.
If 
\beqs
\alpha C_1\big\| \cA^{-1}_1\big\|_{\cH^*\to \cH}\leq 1/2
\eeqs
then $ \matrixA_2^{-1}$ exists,
\begin{align}
&\max\Big\{
\big\| \matrixI - \matrixA_2^{-1}\matrixA_1 \big\|_{\matrixD}
, \big\| \matrixI - \matrixA_1 \matrixA_2^{-1}\big\|_{\matrixD^{-1}}\Big\} \leq 2 \alpha  C_1\big\| \cA^{-1}_1\big\|_{\cH^*\to \cH}
\label{eq:abs1}
\end{align}
and
\beq\label{eq:abs2}
\max\Big\{
\big\| \matrixI - \matrixA_2^{-1}\matrixA_1 \big\|_{2},
\big\| \matrixI - \matrixA_1\matrixA_2^{-1} \big\|_2\Big\} \leq  2\frac{m_+}{m_-} \alpha C_1 \big\| \cA^{-1}_1\big\|_{\cH^*\to \cH}.
\eeq
\end{corollary}


\bre[Implications of Theorem \ref{thm:main} for preconditioning]
If the right-hand side of \eqref{eq:main1} is $\leq c<1$, then 
\bit
\item 
the preconditioned fixed point iteration $\bx^{n+1} = \bx^{n} + \matrixA_2^{-1}({\bf b}-\matrixA_1 \bx^n)$ for solving $\matrixA_1\bx={\bf b}$ satisfies
\beq\label{eq:plant1}
\N{\bx-\bx^n}_{\matrixD} \leq c^n \N{\bx-\bx^0}_{\matrixD},
\eeq
\item when GMRES is applied in the $\matrixD$ inner product, the residual $\br^n:= \matrixA_2^{-1}{\bf b}- \matrixA_2^{-1}\matrixA_1 \bx^n$ satisfies 
\beq\label{eq:plant2}
\frac{\N{\br^n}_{\matrixD}}{\N{\br^0}_{\matrixD}} \leq 
\left(\frac{2 \sqrt{c}}{(1+c)^2}\right)^n
\eeq
(as a consequence of the Elman estimate \cite{El:82,EiElSc:83}; see \cite[Corollary 1.9]{GaGrSp:15}/\cite[Corollary 5.5]{GrPeSp:21}). 
\eit
Similar statements hold for right-preconditioning and/or the fixed-point iteration/GMRES in the Euclidean inner product (using \eqref{eq:main2} instead of \eqref{eq:main1}).

We highlight that, at least in the case of preconditioning by absorption, the action of $\matrixA_2^{-1}$ must be further approximated to provide a practical preconditioner (as discussed in \S\ref{sec:context}).
\ere

\bre[The analogue of Theorem \ref{thm:main} at the level of sesquilinear forms]
A byproduct of the proof of Theorem \ref{thm:main} below is that if
\beq\label{eq:cts1}
 \Big( \N{\mu^{-1}_1-\mu^{-1}_2}_{\cH_0\to \cH_0} 
+ \N{\epsilon_1-\epsilon_2}_{\cH_0\to \cH_0}\Big) \big\| \cA^{-1}_1\big\|_{\cH^*\to \cH}\leq 1/2
\eeq
then $\cA_2^{-1}$ exists and 
\begin{align}\nonumber
&\max\Big\{
\big\| I - \cA_2^{-1}\cA_1 \big\|_{\cH}
, 
\big\| I - \cA_1\cA_2^{-1} \big\|_{\cH}
\Big\}\\
&\hspace{3cm}\leq  2 \Big( \N{\mu^{-1}_1-\mu^{-1}_2}_{\cH_0\to \cH_0} 
+ \N{\epsilon_1-\epsilon_2}_{\cH_0\to \cH_0}\Big) \big\| \cA^{-1}_1\big\|_{\cH^*\to \cH}
\label{eq:cts2}
\end{align}
(see Remark \ref{rem:cts2}).
For the Helmholtz case with $\mu_1=\mu_2$, \cite[Lemma 2.7]{GrPeSp:21} exhibits $\epsilon_1,\epsilon_2$ such that 
the bound \eqref{eq:cts2} is sharp in its $k$-dependence (in this example $\| \cA^{-1}_1\|_{\cH^*\to \cH} \sim k$).
\ere

\bre[Measuring the difference in the coefficients in weaker norms]
By using standard finite-element inverse estimates, \cite[Lemma 5.3]{GrPeSp:21} proved (for the Helmholtz case) bounds on 
$\big\| \matrixI - \matrixA_2^{-1}\matrixA_1 \big\|_{\matrixD}$, 
$\big\| \matrixI - \matrixA_1\matrixA_2^{-1} \big\|_{\matrixD}$, $\big\| \matrixI - \matrixA_1\matrixA_2^{-1} \big\|_{2}$, and 
$\big\| \matrixI - \matrixA_2^{-1}\matrixA_1 \big\|_{2}$ with $L^q$ norms of both $\mu^{-1}_1 - \mu^{-1}_2$ and $\epsilon_1-\epsilon_2$ appearing on the right-hand sides. These bounds also hold for the Helmholtz and Maxwell problems in \S\ref{sec:assumptions}; however, we do not include them here since (i) these bounds involves powers of $h^{-1}$ (from the inverse estimates) and are thus worse than \eqref{eq:main1} and \eqref{eq:main2}, and 
(ii) the bounds 
\eqref{eq:main1} and \eqref{eq:main2} are sufficient for studying preconditioning with absorption (becoming \eqref{eq:abs1} and \eqref{eq:abs2} in Corollary \ref{cor:absorption}) and, in particular, verifying the assumption about this in \cite{LiHu:25}.
\ere

\bre[Theorem \ref{thm:main} improves the Helmholtz results of \cite{GaGrSp:15} and \cite{GrPeSp:21} in three ways]

\red{The bounds \eqref{eq:main1}/\eqref{eq:main2} and \eqref{eq:abs2} are essentially the same as those proved in \cite{GrPeSp:21} and \cite{GaGrSp:15}, respectively. Nevertheless, Theorem \ref{thm:main} improves the results of \cite{GaGrSp:15, GrPeSp:21} by enlarging the class of Helmholtz problems to which the bounds apply, and weakening the assumptions on the finite-dimensional space. In more detail:}

(i) Theorem \ref{thm:main} covers all the Helmholtz problems in \red{Lemma \ref{lem:Helmholtz}}, whereas \cite{GaGrSp:15} and \cite{GrPeSp:21} only consider 
the exterior Dirichlet problem with a nontrapping obstacle (with the obstacle allowed to be empty). For this problem, \cite{GaGrSp:15} considers truncation of the exterior domain by an impedance boundary condition, and \cite{GrPeSp:21} considers truncation by either an impedance boundary condition or the exact Dirichlet-to-Neumann map.

(ii) Theorem \ref{thm:main} makes assumptions only about  $a_1(\cdot,\cdot)$ (i.e., the problem we want to precondition), whereas \cite{GaGrSp:15} and \cite{GrPeSp:21} make assumptions about both $a_1(\cdot,\cdot)$ and $a_2(\cdot,\cdot)$; e.g., the mesh threshold in \cite[Equation 1.12]{GaGrSp:15} is for the problem with absorption, i.e., $a_2(\cdot,\cdot)$, and \cite[Theorems 2.2 and 2.3]{GrPeSp:21} make assumptions on both $a_1(\cdot,\cdot)$ and $a_2(\cdot,\cdot)$.

(iii) Theorem \ref{thm:main} is proved under the discrete inf-sup condition \eqref{eq:dis_both} (which, by Theorem \ref{thm:inverse}, is equivalent to the Galerkin matrix being invertible), whereas the bounds of \cite{GaGrSp:15} and \cite{GrPeSp:21} are proved under stronger assumptions than a discrete inf-sup condition. 
Indeed, 
\cite[Equations 1.13 and 1.14]{GaGrSp:15} prove  bounds analogous to \eqref{eq:main2} under the assumption that, given any data and a sequence of finite-dimensional subspaces, 
the sequence of Galerkin solutions 
is quasi-optimal, with constant independent of $k$ -- we see in Corollary \ref{cor:new1} below that this implies 
that the Galerkin matrix is invertible, with the norm of its inverse 
bounded by the norm of the PDE solution operator.
Furthermore, \cite[Lemma 5.3]{GrPeSp:21} proves bounds analogous to \eqref{eq:main1} and \eqref{eq:main2} under assumptions weaker than quasi-optimality, but stronger than invertibility of the Galerkin matrix; see \cite[Condition 4.2]{GrPeSp:21}.
\label{rem:improve}
\ere

\subsection{\red{Auxiliary result:}~conditions under which the discrete inf-sup constant is proportional to the (continuous) inf-sup constant}\label{sec:is}

The main observation in this subsection (which we couldn't find elsewhere in the literature) is that the argument bounding $\|\mathcal{A}_\ell^{-1}\|_{\cH^*\to \cH}$ in terms of $\|\mathcal{A}_\ell^{-1}\|_{\cH_0\to \cH}$ (in Lemma \ref{lem:is1}) also holds at the discrete level, giving the following result.

\ble[Bound on discrete inf-sup constant via bound on Galerkin solution]\label{lem:new}
Suppose that the assumptions in \S\ref{sec:assumptions} hold and 
there exists $C_{\rm sol,\ell}>0$ such that, for all $f\in \cH_0$, the solution $u_N\in \cH_N$ to 
\beqs
a_\ell(u_N,v_N)= (f,v_N)_{\cH_0}\quad \tfa v_N\in \cH_N
\eeqs
exists and satisfies 
\beq\label{eq:bottle1}
\N{u_N}_{\cH} \leq C_{\rm sol, \ell} \N{f}_{\cH_0}.
\eeq
Then 
\beqs
\inf_{u_{\red N}\in \cH_h\setminus\{0\}} \sup_{v_{\red N}\in \cH_{\red N}\setminus\{0\}} \frac{
\big| a_\ell(u_{\red N},v_{\red N})\big|
}{
\N{u_{\red N}}_{\cH} \N{v_{\red N}}_{\cH}
}
\geq
\frac{1}
{
(C_{\rm G1, \ell})^{-1} \big( 1 + C_{\rm G2, \ell}C_{\rm sol,\ell}\big)}
.
\eeqs
\ele

\red{In other words, the bound \eqref{eq:bottle1} on the discrete solution  for the restricted class of data $f\in \cH_0$ (instead of $f\in \cH^*$) is sufficient to obtain a bound on the discrete inf-sup constant.}

Several analyses of the $h$-FEM for Helmholtz and Maxwell prove the bound \eqref{eq:bottle1} directly; e.g., this is done  in \cite[Theorem 4.1]{LuWu:24}  for the Maxwell impedance problem in \red{Part (ii) of Lemma \ref{lem:Maxwell}}
when $\Omega$ a $C^2$ domain, $\mu_\ell=\epsilon_\ell=I$, $(\cH_N)_{N=0}^\infty$ consists of N\'ed\'elec edge elements of the second type, and $(kh)^2 k$ sufficiently small. 

We now give two other ways of proving the bound \eqref{eq:bottle1}, and hence bounding the discrete inf-sup constant.
The first of these is in the \emph{asymptotic} regime; i.e., when  the sequence of Galerkin solutions is quasi-optimal -- for the $h$-FEM applied to Helmholtz/Maxwell, this is when $(kh)^p \|\cA_\ell^{-1}\|$ is sufficiently small; the second of these is in the \emph{preasymptotic} regime.

\begin{corollary}[Bound on discrete inf-sup constant via quasi-optimality]\label{cor:new1}
Suppose that the assumptions in \S\ref{sec:assumptions} hold and there exists $0<C_{\rm qo,\ell}<\infty$ such that the following holds. Suppose that $(\cH_N)_{N=0}^\infty$ is a sequence of finite-dimensional subspaces such that, for all $N$, given $f\in \cH_0$ 
 the Galerkin solution $u_N$ to 
 \beqs
 a_\ell(u_N,v_N)= (f,v_N)_{\cH_0}\quad\tfa v_N\in \cH_N
 \eeqs
exists, is unique, and satisfies 
 \beqs
 \N{u-u_N}_{\cH}\leq C_{\rm qo,\ell} \N{(I-\Pi_N)u}_{\cH},
 \eeqs
 where $\Pi_N:\cH\to \cH_N$ is the orthogonal projection (in the $\cH$ norm) and $u$ is the solution to $a_\ell(u,v)=(f,v)_{\cH_0}$ for all $v\in \cH$. 

Then 
\beqs
\inf_{u_h\in \cH_h\setminus\{0\}} \sup_{v_h\in \cH_h\setminus\{0\}} \frac{
\big| a_\ell(u_h,v_h)\big|
}{
\N{u_h}_{\cH} \N{v_h}_{\cH}
}
\geq
\frac{1}
{(C_{\rm G1, \ell})^{-1} \Big( 1 + C_{\rm G2, \ell}\big(1 +C_{\rm qo, \ell}\big)\big\|\mathcal{A}^{-1}_\ell\big\|_{\cH_0\to \cH}\Big)}.
\eeqs
\end{corollary}

\begin{corollary}[Bound on discrete inf-sup constant in preasymptotic regime]\label{cor:new2}
Suppose that given $k_0>0$ there exists $c$ such that $\big\|\cA^{-1}_\ell\big\|_{\cH_0\to\cH} \geq c$ for all $k\geq k_0$.
Given $p\in \mathbb{Z}^+$, let $(\cH_h)_{h>0}$ be a sequence of finite-dimensional subspaces of $\cH$.

Suppose there exist $C_j, j=1,2,3,$ such that if 
\beq\label{eq:mesh_threshold}
(kh)^{2p} \big\|\cA^{-1}_\ell\big\|_{\cH_0\to\cH} \leq C_1
\eeq
then, for all $u\in \cH, u_h\in \cH_h$ satisfying $a_\ell(u-u_h, v_h)= 0$ for all $v_h\in \cH_h$,
\beq\label{eq:pre1}
\N{u-u_h}_{\cH} \leq C_2 \Big( 1 + (kh)^p \big\|\cA^{-1}_\ell\big\|_{\cH_0\to\cH}\Big) \N{(I-\Pi_h)u}_\cH,
\eeq
where $\Pi_h:\cH\to \cH_h$ is the orthogonal projection (in the $\cH$ norm), and, 
furthermore
\beq\label{eq:eta}
\big\|(I-\Pi_h)\cA^{-1}_\ell\big\|_{\cH_0\to \cH}\leq C_3 \Big ( kh + (kh)^p \big\| \cA^{-1}_\ell\big\|_{\cH_0\to\cH}\Big).
\eeq

Then, given $0< \varepsilon\leq C_1$ there exists $C_4>0$ such that if $(kh)^{2p} \big\|\cA^{-1}_\ell\big\|_{\cH_0\to\cH}\leq \varepsilon$ then
\beq\label{eq:new2bound}
\inf_{u_h\in \cH_h\setminus\{0\}} \sup_{v_h\in \cH_h\setminus\{0\}} \frac{
\big| a_\ell(u_h,v_h)\big|
}{
\N{u_h}_{\cH} \N{v_h}_{\cH}
}
\geq 
\frac{1}{
(C_{\rm G1, \ell})^{-1} \Big( 1 + C_{\rm G2, \ell} (1+ C_2 C_3 C_4 )\big\|\mathcal{A}^{-1}_\ell\big\|_{\cH_0\to\cH}\Big)
}.
\eeq
\end{corollary}

Since the proofs of Corollaries \ref{cor:new1} and \ref{cor:new2} are so short, we give them here.

\bpf[Proof of Corollary \ref{cor:new1}]
By the triangle inequality and the fact that $\|(I-\Pi_N)u\|_{\cH}\leq \N{u}_\cH$ (since $\Pi_N:\cH\to \cH_N$ is the orthogonal projection),
\beqs
\N{u_N}_{\cH}\leq \N{u}_{\cH} +\N{u-u_N}_{\cH}\leq (1+ C_{\rm qo,\ell}) \N{u}_{\cH}\leq (1+ C_{\rm qo,\ell})\N{\cA^{-1}_\ell}_{\cH_0\to \cH}\N{f}_{\cH_0},
\eeqs
and the result follows from Lemma \ref{lem:new}.
\epf

\bpf[Proof of Corollary \ref{cor:new2}]
Given $f\in \cH_0$, let $u\in \cH$ be the solution to $a_\ell(u,v)=(f,v)_{\cH_0}$ for all $v\in \cH$. Now, by the combination of \eqref{eq:pre1} and \eqref{eq:eta},
\begin{align*}
\N{u-u_h}_{\cH}&\leq C_2 \Big( 1 + (kh)^p \big\|\cA^{-1}_\ell\big\|_{\cH_0\to\cH}\Big)C_3 \Big ( kh + (kh)^p \big\| \cA^{-1}_\ell\big\|_{\cH_0\to\cH}\Big) \N{f}_{\cH_0}.
\end{align*}
Since $(kh)^{2p} \big\| \cA^{-1}_\ell\big\|_{\cH_0\to\cH}\leq \varepsilon$
and $\big\| \cA^{-1}_\ell\big\|_{\cH_0\to\cH}\geq c$, there exists $C_4>0$ (depending on $c$ and $\varepsilon$) such that 
\beqs
\N{u-u_h}_{\cH}\leq C_2 C_3 C_4 \big\| \cA^{-1}_\ell\big\|_{\cH_0\to\cH}\N{f}_{\cH_0}.
\eeqs
By the triangle inequality, \eqref{eq:bottle1} holds with 
\beqs
C_{\rm sol, \ell}:= (1 + C_2 C_3 C_4)\big\| \cA^{-1}_\ell\big\|_{\cH_0\to\cH},
\eeqs
and the result \eqref{eq:new2bound} follows from Lemma \ref{lem:new}.
\epf

\bre[$\cA^{-1}_\ell$ vs $(\cA^*_\ell)^{-1}$ in the bound \eqref{eq:eta}]
The bound \eqref{eq:eta} is usually proved with $\cA^{-1}_\ell$ replaced by $(\cA^*_\ell)^{-1}$ (see, e.g., \cite{MeSa:10, ChNi:20, GS3}). However, if 
\beq\label{eq:plane1}
a(u,v)= a(\overline{v},\overline{u}) \quad\tfa u,v \in \cH,
\eeq
then the definitions of $\cA^*_\ell$ and $\cA_\ell$ imply that $(\cA^*_\ell)^{-1} f = \overline{\cA_\ell^{-1} \overline{f}}$, and thus the bound \eqref{eq:eta} is equivalent to the corresponding bound with $\cA^{-1}_\ell$ replaced by $(\cA^*_\ell)^{-1}$. The condition \eqref{eq:plane1} holds 
for all the Helmholtz and Maxwell problems in \red{Lemmas \ref{lem:Helmholtz} and \ref{lem:Maxwell}} \emph{either} 
when $\mu_\ell^{-1}$ and $\epsilon_\ell$ are Hermitian \emph{or} when $\mu_\ell^{-1}$ and $\epsilon_\ell$ corresponding to a radial PML (see \cite[Lemma 2.3]{GLSW1}).

Furthermore, the proof of \cite[Theorem 1.7]{GS3}, which establishes \eqref{eq:eta} with $\cA^{-1}_\ell$ replaced by $(\cA^*_\ell)^{-1}$ for general Helmholtz problems, also proves \eqref{eq:eta}.
\ere

\bre[Link to other results in the literature]\label{rem:link}
The result that the discrete inf-sup constant is bounded by the continuous inf-sup constant when $\|(I-\Pi_h)(\cA^*_\ell)^{-1}\|_{\cH_0\to \cH}$ is small is \cite[Theorem 4.2]{MeSa:10}. By the Schatz argument \cite{Sc:74}, \cite[Theorem 4.3]{MeSa:10}, if $\|(I-\Pi_h)(\cA^*_\ell)^{-1}\|_{\cH_0\to \cH}$ is small, then the sequence of Galerkin solutions is quasi-optimal, and so  \cite[Theorem 4.2]{MeSa:10} is morally equivalent to Corollary \ref{cor:new1}.
\ere

\paragraph{Bounding the discrete inf-sup constant for the $h$-FEM applied to each of the Helmholtz and Maxwell problems in \S\ref{sec:assumptions}.}

Suppose that the coefficients and domain satisfy the natural regularity requirements for using degree $p$ polynomials; i.e., the domain is $C^{p,1}$ and the coefficients are piecewise $C^{p-1,1}$.
\bit
\item All the Helmholtz problems in \red{Lemma \ref{lem:Helmholtz}} satisfy the assumptions of Corollary \ref{cor:new2} by \cite[Theorems 1.5 and 1.7]{GS3} (with \cite[Theorem 1.5]{GS3} giving \eqref{eq:pre1} and \cite[Theorem 1.7]{GS3} giving \eqref{eq:eta}).
\item For the Maxwell PEC problem of \red{Part (i) of Lemma \ref{lem:Maxwell}} discretised using N\'ed\'elec edge elements of either first or second type, the bound \eqref{eq:pre1} under the condition \eqref{eq:mesh_threshold} is proved in \cite[Theorem 1.3]{CGS1}. 
The bound \eqref{eq:eta} is not proved in \cite{CGS1} (the arguments of \cite{CGS1} use a duality argument that ``builds in" the splitting used to prove \eqref{eq:eta}). Therefore, although we expect that the 
bound \eqref{eq:bottle1} holds in the preasymptotic regime for the Maxwell PEC problem (via a result analogous to Corollary \ref{cor:new2}), right now \cite[Theorem 1.3]{CGS1} only allows us to use Corollary \ref{cor:new1} in the asymptotic regime (i.e.~under the condition that $(kh)^p \|\cA^{-1}_\ell\|$ is sufficiently small). 
\item As noted above, for the Maxwell impedance problem of \red{Part (ii) of Lemma \ref{lem:Maxwell}}, the bound \eqref{eq:bottle1} is proved directly in \cite[Theorem 4.1]{LuWu:24}  for the Maxwell impedance problem in \red{Part (ii) of Lemma \ref{lem:Maxwell}} 
when $\Omega$ a $C^2$ domain, $\mu_\ell=\epsilon_\ell=I$, $(\cH_N)_{N=0}^\infty$ consists of N\'ed\'elec edge elements of the second type, and $(kh)^2 k$ sufficiently small (note that $\|\cA^{-1}_\ell\|\sim k$ here by \cite{HiMoPe:11a}).
\eit

\section{Proofs of the main results}\label{sec:proofs}

\red{Sections \ref{sec:5.1}-\ref{sec:5.3} are devoted to the proof of Theorem \ref{thm:main}, and Section \ref{sec:Geneva} is devoted to the proof of Lemma \ref{lem:new}.} 

\subsection{Bounds on $\| \matrixI - \matrixA_2^{-1}\matrixA_1\|$
and $\| \matrixI -\matrixA_1 \matrixA_2^{-1}\|$ in terms of the discrete inf-sup constant
}\label{sec:5.1}

The heart of the proof of Theorem \ref{thm:main} is the following result.

\begin{lemma}
\label{lem:main}
Suppose that the assumptions of \S\ref{sec:assumptions} hold. 
Given $N>0$, if  exists $0<\Cdiso<\infty$ such that
\begin{subequations}\label{eq:disold_both}
\beq\label{eq:disold}
\inf_{u_N\in \cH_N\setminus\{0\}} \sup_{v_N\in \cH_N\setminus\{0\}} \frac{
\big| a_2(u_N,v_N)\big|
}{
\N{u_N}_{\cH} \N{v_H}_{\cH}
}
\geq 
\frac{1}{\Cdiso}
\eeq
and
\beq\label{eq:distold}
\tfa v_N\in \cH_N\setminus \{0\},\, \sup_{u_N\in \cH_N\setminus\{0\}} |a_2(u_N,v_N)|>0,
\eeq
\end{subequations}
then $ \matrixA_2^{-1}$ exists and 
\beq\label{eq:main1old}
\max\Big\{
\big\| \matrixI - \matrixA_2^{-1}\matrixA_1 \big\|_{\matrixD}
, \big\| \matrixI - \matrixA_1 \matrixA_2^{-1}\big\|_{\matrixD^{-1}}\Big\}
\leq  \Big( \N{\mu^{-1}_1-\mu^{-1}_2}_{\cH_0\to \cH_0} 
+ \N{\epsilon_1-\epsilon_2}_{\cH_0\to \cH_0}\Big) \Cdiso.
\eeq
Furthermore, if $\mu_1=\mu_2$, then 
\beq\label{eq:main2old}
\max\Big\{
\big\| \matrixI - \matrixA_2^{-1}\matrixA_1 \big\|_{2},
\big\| \matrixI - \matrixA_1\matrixA_2^{-1} \big\|_2\Big\} \leq \frac{m_+}{m_-}\N{\epsilon_1-\epsilon_2}_{\cH_0\to \cH_0} \Cdiso.
\eeq
\end{lemma}
 
Given Lemma \ref{lem:main}, to prove Theorem \ref{thm:main} we only need to 
show that the condition \eqref{eq:dis_both} implies that \eqref{eq:disold_both} holds with $\Cdiso$ bounded above by $\|\cA^{-1}_1\|_{\cH^*\to \cH}$. 
 
 \bpf[Proof of Lemma \ref{lem:main}]
We first show how the bounds on $\matrixI - \matrixA_1\matrixA_2^{-1}$ follow from those on $\matrixI - \matrixA_2^{-1}\matrixA_1$.
Given a matrix $\matrixC$, let $\matrixC^\dagger$ be the conjugate transpose of $\matrixC$ (i.e.~the adjoint with respect to $(\cdot,\cdot)_2$). 
Then
\beqs
 \big\| \matrixI - \matrixA_1 \matrixA_2^{-1}\big\|_{2}
 =
\big\| \matrixI - (\matrixA_2^\dagger)^{-1}(\matrixA_1^\dagger) \big\|_{2}.
\eeqs
Furthermore, by direct calculation,
\beqs
\frac{
(\bV_1, \matrixC \bV_2)_{\matrixD^{-1}}
}{
\N{\bV_1}_{\matrixD^{-1}}\N{\bV_2}_{\matrixD^{-1}}
} 
=\frac{
(\matrixC^\dagger \bW_1, \bW_2)_{\matrixD}
}{
\N{\bW_1}_{\matrixD}\N{\bW_2}_{\matrixD}
}  
\quad\tfa \bV_j\in \Com^N,
\eeqs
where $\bW_j := \matrixD^{-1}\bV_j$, $j=1,2,$.
Therefore
\beqs
 \big\| \matrixI - \matrixA_1 \matrixA_2^{-1}\big\|_{\matrixD^{-1}}
 =
\big\| \matrixI - (\matrixA_2^\dagger)^{-1}(\matrixA_1^\dagger) \big\|_{\matrixD}.
\eeqs
Recall from \S\ref{sec:assumptions} that if $a_\ell(\cdot,\cdot)$ satisfies the assumptions in \S\ref{sec:assumptions}, then so does $a^*_\ell(\cdot,\cdot)$ defined by $a^*_\ell(u,v)=
\overline{a_\ell(v,u)}$, with the same constants $C_{\rm G1,\ell}$ and $C_{\rm G2,\ell}$. 
Furthermore, if the discrete inf-sup condition \eqref{eq:dis_both}  on $a(\cdot,\cdot)$ holds, then it also holds for $a^*(\cdot,\cdot)$, with the same $\Cdiso$; see, e.g., \cite[Lemma 6.5.3]{Ha:92}, \cite[Remark 2.1.45]{SaSc:11} (this is just the result that the norm of the adjoint operator equals the norm of the original operator).

Therefore, once the bound on $\big\| \matrixI - \matrixA_2^{-1}\matrixA_1 \big\|_{\matrixD}$ in \eqref{eq:main1} is established, this bound also holds $\matrixA_2$ replaced by $\matrixA_2^\dagger$ and $\matrixA_1$ replaced by $\matrixA_1^\dagger$, so that 
\begin{align*}
 \big\| \matrixI - \matrixA_1 \matrixA_2^{-1}\big\|_{\matrixD^{-1}}
& =\big\| \matrixI - (\matrixA_2^\dagger)^{-1}(\matrixA_1^\dagger) \big\|_{\matrixD}\\
 &\leq  \Big( \N{(\mu^*_1)^{-1}-(\mu^*_2)^{-1}}_{\cH_0\to \cH_0} 
+ \N{\epsilon^*_1-\epsilon^*_2}_{\cH_0\to \cH_0}\Big) \Cdiso\\
&=  \Big( \N{\mu^{-1}_1-\mu^{-1}_2}_{\cH_0\to \cH_0} 
+ \N{\epsilon_1-\epsilon_2}_{\cH_0\to \cH_0}\Big) \Cdiso.
\end{align*}
Identical reasoning obtains the bound on $\| \matrixI - \matrixA_1 \matrixA_2^{-1}\|_{2}$ from that on $\| \matrixI -  \matrixA_2^{-1}\matrixA_1\|_{2}$.

We now prove the bounds on $\matrixI - \matrixA_2^{-1}\matrixA_1$ in \red{\eqref{eq:main1old}} and \red{\eqref{eq:main2old}}. 
By Theorem \ref{thm:inverse}, the discrete inf-sup condition \eqref{eq:dis_both} implies that $\matrixA_2: \Com^{N}\to \Com^N$ is invertible. 
Then, by the definitions of $\matrixA_\ell$ \eqref{eq:matrixA}, $\matrixS_\mu$, and $\matrixM_\epsilon$ \eqref{eq:matrixS},
\beqs
\matrixI - \matrixA_2^{-1}\matrixA_1 = \matrixA_2^{-1} (\matrixA_2 -\matrixA_1) = \matrixA_2^{-1} \Big( \matrixS_{\mu^{-1}_2} -\matrixS_{\mu^{-1}_1} - \matrixM_{\epsilon_2}+ \matrixM_{\epsilon_1}\Big)
= \matrixA_2^{-1} \Big( \matrixS_{\mu^{-1}_2-\mu^{-1}_1} - \matrixM_{\epsilon_2-\epsilon_1}\Big).
\eeqs
Therefore, to prove \eqref{eq:main1old}, it is sufficient to prove that 
\beqs
\big\| \matrixA_2^{-1} \matrixS_{\mu^{-1}} \big\|_{\matrixD} \leq \N{\mu^{-1}}_{\cH_0\to \cH_0}\Cdiso \quad\tand\quad
\big\| \matrixA_2^{-1} \matrixM_\epsilon \big\|_{\matrixD} \leq \N{\epsilon}_{\cH_0\to \cH_0}\Cdiso\red{,}
\eeqs
and to prove \eqref{eq:main2old}, it is sufficient to prove that 
\beqs
\big\| \matrixA_2^{-1} \matrixM_\epsilon \big\|_{2} \leq \frac{m_+}{m_-}\N{\epsilon}_{\cH_0\to \cH_0}\Cdiso.
\eeqs
It is therefore sufficient to prove that, given $\bF\in \Com^N$, the solutions $\bU$ and $\bW$ to 
\beq\label{eq:last1}
\matrixA_2 \bU = \matrixS_{\mu^{-1}} \bF\quad\tand\quad\matrixA_2\bW=\matrixM_\epsilon \bF
\eeq
satisfy
\beq\label{eq:asleep1}
\N{\bU}_{\matrixD}\leq \N{\mu^{-1}}_{\cH_0\to \cH_0} \Cdiso \N{\bF}_{\matrixD},\quad
\N{\bW}_{\matrixD} 
 \leq \N{\epsilon}_{\cH_0\to \cH_0}\Cdiso\N{\bF}_{\matrixD},
\eeq
and 
\beq\label{eq:asleep2}
m_-\N{\bW}_{2} 
 \leq m_+\N{\epsilon}_{\cH_0\to \cH_0}\Cdiso\N{\bF}_{2}.
\eeq
Given $\bF\in \Com^N$, let $\widetilde{f} \in \cH$ and $\widetilde{F} \in \cH^*$ be defined by
\beq\label{eq:specialdata}
\widetilde{f} := \sum_{j=1}^N F_j \phi_j
\quad\tand\quad 
\widetilde{F}(v) = \big( \mu^{-1} \cD \widetilde{f} ,\cD v\big)_{\cH_0}.
\eeq
The solutions $u_N$ and $w_N$ to the variational problems
\beq\label{eq:uNprob}
a_2(u_N,v_N)= \widetilde{F}(v_N) \quad\tand\quad a_2(w_N,v_N) = (\widetilde{f}, v_N)_{\cH_0} \quad\tfa v_N \in \cH_N
\eeq
are then such that 
\beq\label{eq:final1}
u_N= \sum_{j=1}^N U_j \phi_j \quad\tand\quad w_N= \sum_{j=1}^N W_j \phi_j,
\eeq
with $\bU$ and $\bW$ the solutions to \eqref{eq:last1}.

Now, by the bound \eqref{eq:disold} involving $\Cdiso$,  the definition of $u_N$ \eqref{eq:uNprob}, the definition of $\widetilde{F}$ \eqref{eq:specialdata}, and the fact that $\N{\cD}_{\cH\to \cH_0}\leq 1$, 
\begin{align*}
\N{u_N}_{\cH} \leq \Cdiso \sup_{v_N\in \cH_N} \frac{
\big|a_2(u_N,v_N)\big|
}{
\N{v_N}_{\cH}
}
&\leq 
\Cdiso \sup_{v_N\in \cH_N} \frac{ \big\|\mu^{-1} \cD \widetilde{f}\big\|_{\cH_0}\N{\cD v_N}_{\cH_0}}
{\N{v_N}_{\cH}}\\
&\leq 
\Cdiso\N{\mu^{-1}}_{\cH_0\to \cH_0}
\big\|\widetilde{f}\big\|_{\cH}.
\end{align*}
By \eqref{eq:matrixD}, \eqref{eq:final1}, and \eqref{eq:specialdata}, $\N{u_N}_{\cH}= \N{\bU}_{\matrixD}$ and $\big\|\widetilde{f}\big\|_{\cH}= \N{\bF}_{\matrixD}$, and the first bound in \eqref{eq:asleep1} follows.
Similarly, by the bound \eqref{eq:disold} involving $\Cdiso$, the definitions of $w_N$ \eqref{eq:uNprob}, 
 and the fact that $\N{\cdot}_{\cH_0}\leq \N{\cdot}_\cH$, 
\begin{align}\nonumber
\N{w_N}_{\cH} \leq \Cdiso \sup_{v_N\in \cH_N} \frac{
\big|a_2(w_N,v_N)\big|
}{
\N{v_N}_{\cH}
}
&\leq 
\Cdiso \sup_{v_N\in \cH_N} \frac{ \big\|\epsilon \widetilde{f}\big\|_{\cH_0}\N{ v}_{\cH_0}}
{\N{v}_{\cH}}\\
&\leq 
\Cdiso\N{\epsilon}_{\cH_0\to \cH_0} \big\|\widetilde{f}\big\|_{\cH_0}.
\label{eq:asleep3}
\end{align}
The second bound in \eqref{eq:asleep1} then follows since, by \eqref{eq:matrixD}, \eqref{eq:final1}, and \eqref{eq:specialdata},
 $\N{w_N}_{\cH}= \N{\bW}_{\matrixD}$ and $\|\widetilde{f}\|_{\cH_0}\leq\|\widetilde{f}\|_{\cH}= \N{\bF}_{\matrixD}$. 

Finally, the bound \eqref{eq:asleep2} follows from \eqref{eq:asleep3} by using $m_-\|\bW\|_{2}\leq \|w_N\|_{\cH_0}\leq \|w_N\|_{\cH}$ and $\|\widetilde{f}\|_{\cH_0}\leq m_+ \| \bF\|_2$ (with both of these bounds following from \eqref{eq:normequiv}).
\epf 
 
  \subsection{Norm of solution operator under perturbation}\label{sec:perturb}

\ble[Norm of solution operator under perturbation]\label{lem:perturb}
Suppose that the assumptions of \S\ref{sec:assumptions} hold. 
If $\cA_1:\cH\to \cH^*$ is invertible and 
\beq\label{eq:beer3}
\Big(\N{\mu^{-1}_1-\mu^{-1}_2}_{\cH_0\to \cH_0} 
 + \N{\epsilon_1-\epsilon_2}_{\cH_0\to \cH_0}\Big) \|\mathcal{A}_1^{-1}\|_{\cH^*\to \cH}\leq 1/2,
 \eeq
 then 
 \beq\label{eq:beer3a}
 \|\mathcal{A}_2^{-1}\|_{\cH^*\to \cH}\leq 2 \|\mathcal{A}_1^{-1}\|_{\cH^*\to \cH}.
 \eeq
\ele 

\bpf
Given $F\in \cH^*$, let $u_2\in \cH$ be the solution to $a_2(u_2,v)=F(v)$ for all $v\in \cH$. By the definition of $a_\ell(\cdot,\cdot)$ \eqref{eq:sesqui}, 
\beqs
a_1(u_2, v) = F(v) + \big((\mu^{-1}_1-\mu^{-1}_2) \cD u_2, \cD v\big)_{\cH_0} - \big( (\epsilon_1-\epsilon_2) u_2, v\big)_{\cH_0}.
\eeqs
By the definition of $\|\mathcal{A}_1^{-1}\|_{\cH^*\to \cH}$, the fact that $\|\cD\|_{\cH\to \cH_0}\leq 1$, and the bound \eqref{eq:beer3},
\begin{align*}
\N{u_2}_{\cH} &\leq 
\|\mathcal{A}_1^{-1}\|_{\cH^*\to \cH}\Big(
\N{F}_{\cH^*} + 
\N{\mu^{-1}_1-\mu^{-1}_2}_{\cH_0\to \cH_0} 
\N{u_2}_{\cH} + \N{\epsilon_1-\epsilon_2}_{\cH_0\to \cH_0}\N{u_2}_{\cH_0}\Big)\\
&\leq \|\mathcal{A}_1^{-1}\|_{\cH^*\to \cH}
\N{F}_{\cH^*} + 
\frac{1}{2}\N{u_2}_{\cH}, 
\end{align*}
and the result \eqref{eq:beer3a} follows.
\epf

\bre[Proof of the bound \eqref{eq:cts2} under the condition \eqref{eq:cts1}]\label{rem:cts2}
The condition \eqref{eq:cts1} is the same as \eqref{eq:beer3}, and thus \eqref{eq:beer3a} holds. Since
\beqs
I - \cA_2^{-1}\cA_1 =  \cA_2^{-1}(\cA_2- \cA_1) \quad\tand\quad I - \cA_1\cA_2^{-1} =  (\cA_2- \cA_1)\cA_2^{-1}.
\eeqs
the bound \eqref{eq:cts2} follows from \eqref{eq:beer3a} and the bound 
\beqs
\N{\cA_2- \cA_1}_{\cH\to \cH^*}\leq 
\N{\mu^{-1}_1-\mu^{-1}_2}_{\cH_0\to \cH_0} 
 + \N{\epsilon_1-\epsilon_2}_{\cH_0\to \cH_0};
 \eeqs
 this last bound holds since, by the definition of $a_\ell(\cdot,\cdot)$ \eqref{eq:sesqui}, 
\beqs
\big| a_1(u,v) - a_2(u,v)\big| \leq \Big(\N{\mu^{-1}_1-\mu^{-1}_2}_{\cH_0\to \cH_0} 
 + \N{\epsilon_1-\epsilon_2}_{\cH_0\to \cH_0}\Big) \N{u}_\cH \N{v}_{\cH}
\eeqs
for all $u,v\in \cH$ (since $\| \cD \|_{\cH \to \cH_0}\leq 1$).
\ere

\ble[Norm of discrete solution operator under perturbation]\label{lem:perturb2}
Suppose that the assumptions of \S\ref{sec:assumptions} hold. 
Suppose that, with $0<\Cdisoo<\infty$,
\begin{subequations}\label{eq:beer_both}
\beq\label{eq:beer1}
\inf_{u_N\in \cH_N\setminus\{0\}} \sup_{v_N\in \cH_N\setminus\{0\}} \frac{
\big| a_1(u_N,v_N)\big|
}{
\N{u_N}_{\cH} \N{v_H}_{\cH}
}
\geq 
\frac{1}{\Cdisoo},
\eeq
\beq\label{eq:beer2}
\tfa v_N\in \cH_N\setminus \{0\},\, \sup_{u_N\in \cH_N\setminus\{0\}} |a_1(u_N,v_N)|>0,
\eeq
\end{subequations}
and 
\beq\label{eq:beer4}
\Big(\N{\mu^{-1}_1-\mu^{-1}_2}_{\cH_0\to \cH_0} 
 + \N{\epsilon_1-\epsilon_2}_{\cH_0\to \cH_0}\Big)\Cdisoo\leq 1/2.
 \eeq
Then \eqref{eq:disold} and \eqref{eq:distold} hold with 
 \beqs
\Cdiso=2 \Cdisoo.
 \eeqs
\ele 

\bpf
Let $\cA_{N,\ell}:\cH_N\to (\cH_N)^*$, $\ell=1,2,$ be the operators associated to the sesquilinear forms $a_\ell(\cdot,\cdot):\cH_N\times \cH_N\to \Com$, $\ell=1,2$. 
By Theorem \ref{thm:inverse}, the condition \eqref{eq:beer_both} is equivalent to the statement that 
\beq\label{eq:beer5}
\|\cA_{N,1}^{-1}\|_{(\cH_N)^*\to \cH_N}\leq \Cdisoo.
\eeq
The condition \eqref{eq:beer4} then implies that \eqref{eq:beer3} holds with $\|\mathcal{A}_1^{-1}\|_{\cH^*\to \cH}$ replaced by $\|\cA_{N,1}^{-1}\|_{(\cH_N)^*\to \cH_N}$. 
The arguments in the proof of Lemma \ref{lem:perturb} then show that 
\beqs
\big\|\cA_{N,2}^{-1}\big\|_{(\cH_N)^*\to \cH_N}\leq 2\big\|\cA_{N,1}^{-1}\big\|_{(\cH_N)^*\to \cH_N}.
\eeqs
Then, by Theorem \ref{thm:inverse} applied to $\cA_{N,2}$ and \eqref{eq:beer5},
\beqs
\inf_{u_N\in \cH_N\setminus\{0\}} \sup_{v_N\in \cH_N\setminus\{0\}} \frac{
\big| a_2(u_N,v_N)\big|
}{
\N{u_N}_{\cH} \N{v_H}_{\cH}
}
\geq 
\frac{1}{2\big\|\cA_{N,1}^{-1}\big\|_{(\cH_N)^*\to \cH_N}} \geq \frac{1}{2\Cdisoo},
\eeqs
and the result follows.
\epf
 
 \subsection{Proof of Theorem \ref{thm:main}} \label{sec:5.3}
  
The assumption \eqref{eq:dis_both} implies that \eqref{eq:beer_both} holds with $\Cdisoo= C_1 \|\cA^{-1}_1\|_{\cH^*\to \cH}$.
The condition \eqref{eq:FridaySun1} is then \eqref{eq:beer4} and thus 
Lemma \ref{lem:perturb2} implies that \eqref{eq:disold} and \eqref{eq:distold} hold with 
 \beqs
\Cdiso=2 C_1 \|\cA^{-1}_1\|_{\cH^*\to \cH}.
 \eeqs
 With this value of $\Cdiso$, the bounds \eqref{eq:main1old} and \eqref{eq:main2old} from Lemma \ref{lem:main} then become the results \eqref{eq:main1} and \eqref{eq:main2}.
  
 \subsection{Proofs of Lemma \ref{lem:is1} and \ref{lem:new} (about equivalence of the norms of the continuous and discrete inverses)}\label{sec:Geneva}

We first prove Lemma \ref{lem:is1}, which consists of proving the bounds \eqref{eq:sunny1} and \eqref{eq:sunny2}. 

\paragraph{Proof of \eqref{eq:sunny1}.}
We highlight that this argument goes back to at least \cite[Lemma 3.4]{ChMo:08}, but since it is short, and crucial to the proof of Lemma \ref{lem:new}, we include it for completeness.
The first inequality in \eqref{eq:sunny1} is an immediate consequence of the fact that $\|\cdot\|_{\cH^*}\leq \|\cdot\|_{\cH_0}$ (which follows from the definition of $\|\cdot\|_{\cH^*}$ and $\|\cdot\|_{\cH_0}\leq \|\cdot\|_{\cH}$).
Let $a^+_\ell(u,v):= a_\ell(u,v)+ C_{\rm G2,\ell}(u,v)_{\cH_0}$ and observe that $a^+_\ell(\cdot,\cdot)$ is coercive by \eqref{eq:Ccoer}.
To prove the second inequality in \eqref{eq:sunny1} it is sufficient to prove that, given $F\in \cH^*$, the solution of $a_\ell(u,v)=F(v)$ for all $v\in \cH$ satisfies
\beq\label{eq:STP1}
\N{u}_\cH \leq (C_{\rm G1})^{-1}\big(1+ C_{\rm G2,\ell}\big\|\cA^{-1}_\ell\big\|_{\cH_0\to \cH}\big)\N{F}_{\cH^*}.
\eeq
Let $u^\pm$ by the solutions to 
\beqs
a^+_\ell(u^+,v) =F(v) \quad\tand\quad 
a_\ell(u^-,v)= C_{\rm G2}(u^+,v)_{\cH_0}\quad\tfa v\in \cH;
\eeqs
these definitions imply that $u=u^++u^-$. By \eqref{eq:Ccoer} and the Lax--Milgram lemma,
\beq\label{eq:air1}
\big\|u^+\big\|_{\cH} \leq (C_{\rm G1})^{-1}\N{F}_{\cH^*}.
\eeq
By the definition of $\|\cA_\ell^{-1}\|_{\cH_0\to \cH}$ and the bound \eqref{eq:air1} on $u^+$,
\beq\label{eq:toolate1}
\big\|u^-\big\|_{\cH} \leq \|\cA_\ell^{-1}\|_{\cH_0\to \cH} C_{\rm G2,\ell} \big\|u^+\big\|_{\cH_0} \leq \|\cA_\ell^{-1}\|_{\cH_0\to \cH} C_{\rm G2,\ell} (C_{\rm G1})^{-1}\N{F}_{\cH^*};
\eeq
the bound \eqref{eq:STP1} -- and hence also the second bound in \eqref{eq:sunny1} -- then follows.

\paragraph{Proof of \eqref{eq:sunny2}.}
The first inequality in \eqref{eq:sunny2} is an immediate consequence of the fact that $\|\cdot\|_{\cH_0}\leq \|\cdot\|_{\cH}$. To prove the second inequality, observe that the G\aa rding-type inequality \eqref{eq:Ccoer} implies that, given $f\in \cH_0$, 
\begin{align*}
C_{\rm G1,\ell} \big\| \mathcal{A}^{-1}_\ell f \big\|_{\cH}^2 &\leq C_{\rm G2, \ell} \big\| \mathcal{A}^{-1}_\ell f\big\|^2_{\cH_0} + \big| \langle \mathcal{A}^{-1}_\ell f, f \rangle \big| ,\\
&\leq C_{\rm G2, \ell} \big\| \mathcal{A}^{-1}_\ell f\big\|^2_{\cH_0} +\big\|\mathcal{A}^{-1}_\ell f\big\|_{\cH_0} \|f \|_{\cH_0}.
\end{align*}
Therefore, 
\beqs
C_{\rm G1,\ell} \big\| \mathcal{A}^{-1}_\ell \big\|_{\cH_0\to \cH}^2 
\leq C_{\rm G2, \ell} \big\| \mathcal{A}^{-1}_\ell \big\|^2_{\cH_0\to \cH_0} +\big\|\mathcal{A}^{-1}_\ell \big\|_{\cH_0\to\cH_0},
\eeqs
and the second bound in \eqref{eq:sunny2} follows.

\

Finally, to prove Lemma \ref{lem:new} we observe the argument above proving the second bound in \eqref{eq:sunny1} remains unchanged when $F\in \cH^*$ is replaced by $F\in (\cH_N)^*$; indeed, the Lax--Milgram lemma proves that
\beqs
\big\|u^+_N\big\|_{\cH} \leq (C_{\rm G1,\ell})^{-1}\N{F}_{(\cH_N)^*}.
\eeqs
Furthermore, since $u_N^+ \subset \cH_N\subset \cH\subset\cH_0$, \eqref{eq:bottle1} implies that 
\beqs
\big\|u^-_N\big\|_{\cH} \leq C_{\rm sol,\ell} C_{\rm G2,\ell} \big\|u^+\big\|_{\cH_0} \leq C_{\rm sol,\ell} C_{\rm G2,\ell}(C_{\rm G1,\ell})^{-1}\N{F}_{(\cH_N)^*}
\eeqs
(compare to \eqref{eq:toolate1}). The result of Lemma \ref{lem:new} therefore follows.

\section*{Acknowledgements}

The author thanks Qiya Hu (Academy of Mathematics and Systems Science, Chinese Academy of Sciences) for asking him the question of whether the results of \cite{GaGrSp:15, GrPeSp:21} extend to the time-harmonic Maxwell equations, Théophile Chaumont-Frelet (INRIA, Lille) for discussions about the discrete inf-sup constant of the Helmholtz sesquilinear form, \red{and the referees for many useful comments.}

\footnotesize{
\bibliographystyle{plain}
\bibliography{combined.bib}
}

\end{document}